\documentclass[10pt,A4,reqno]{amsart}
\usepackage{amsmath,amstext,amsbsy,amsopn,amsthm,%
amscd,amsfonts,amssymb}
\usepackage{mathrsfs}

\numberwithin{equation}{section}

\theoremstyle{definition} %
\newtheorem{defn}{Definition}[section]
\newtheorem{prop}[defn]{Proposition}
\newtheorem{thm}[defn]{Theorem}
\theoremstyle{remark}
\newtheorem{remk}[defn]{Remark}

\newtheorem{corl}[defn]{Corollary}
\newtheorem*{ack}{Acknowledgements}

\newcommand{\REMOVE}[1]{}

\newcommand{\changed}[1]{#1}
\newcommand{\changedd}[1]{#1}

\newcommand{\eps}{\varepsilon}

\newcommand{\gO}{\varOmega}
\newcommand{\ofO}{{(\gO)}}
\newcommand{\ba}{\begin{align*}}
\newcommand{\ea}{{\end{align*}}}
\newcommand{\les}{\leqslant}
\newcommand{\ges}{\geqslant}
\newcommand{\tto}{\longrightarrow}
\newcommand{\abs}[1]{\lvert#1\rvert}
\newcommand{\comp}{\raisebox{0.75pt}{\mbox{$\:{\scriptstyle \circ}\:$}}}
\newcommand{\cl}[1]{{\overline{#1}}}
\newcommand{\absv}[1]{{\langle#1\rangle}}
\newcommand{\lin}{T}

\newcommand{\cee}[1]{C^{#1}\!}
\newcommand{\lone}{{L^1}}
\newcommand{\linf}{{L^\infty}}
\newcommand{\sob}{{W^1_1}}
\newcommand{\rigs}{\mathscr{R}}
\newcommand{\qrigs}{\!/\!\rigs}
\newcommand{\ld}{{LD}}
\newcommand{\ldo}{\ld\ofO}
\newcommand{\trace}{\gamma}

\newcommand{\bndo}{{\partial\gO}}
\newcommand{\cloo}{\cl{\gO}}
\newcommand{\dis}{\chi}
\newcommand{\sts}{\varSigma}
\newcommand{\ebdfc}{T}
\newcommand{\ofbdo}{(\bndo)}
\newcommand{\ofclo}{(\cloo)}
\newcommand{\vono}{(\gO,\rthree)}
\newcommand{\vonbdo}{(\bndo,\rthree)}
\newcommand{\vonclo}{(\cl{\gO},\rthree)}
\newcommand{\strono}{(\gO,\reals^{6})}

\newcommand{\infclo}{{\infty,\cloo}}
\newcommand{\infbdo}{{\infty,\bndo}}
\newcommand{\onebdo}{{1,\bndo}}
\newcommand{\pr}{\proj_{\rigs}}
\newcommand{\pq}{\proj}
\newcommand{\qr}{\!/\!\reals}
\newcommand{\nvf}{\mathfrak{n}}
\newcommand{\nvfs}{{\mathfrak{N}}}
\newcommand{\mnorm}{\abs}
\newcommand{\rig}{r}

\newcommand{\reals}{\ensuremath{\mathbb{R}}}
\newcommand{\rthree}{\ensuremath{\reals^{\scriptscriptstyle 3}}}

\DeclareMathOperator{\ess}{ess}
\DeclareMathOperator*{\esssup}{ess\;sup}
\DeclareMathOperator{\kernel}{Kernel}
\newcommand{\from}{\colon}
\DeclareMathOperator{\image}{Image}
\newcommand{\proj}{\pi}
\newcommand{\vs}{\ensuremath{\mathbf{W}}}

\newcommand{\ino}[1]{\int\limits_{#1}}
\newcommand{\half}{\ensuremath{{\frac{1}{2}}}}

\newcommand{\paren}[1]{\left(#1\right)}

\newcommand{\braces}[1]{\left\{#1\right\}}
\newcommand{\norm}[1]{\lVert#1\rVert}
\newcommand{\fun}{\varphi}

\newcommand{\pis}{\ensuremath{x}}

\newcommand{\vect}{\ensuremath{v}}
\newcommand{\vbase}{e}

\newcommand{\bdry}{\partial}
\newcommand{\nor}{\ensuremath{{\nu}}}

\newcommand{\dV}{}
\newcommand{\dA}{}

\newcommand{\vf}{w}

\newcommand{\fc}{F}
\newcommand{\bfc}{{b}}
\newcommand{\sfc}{{t}}
\newcommand{\st}{\sigma}
\begin{document}
\title{Bounds on the Trace Mapping of
 LD-Fields}

\author{Ronen Peretz and Reuven Segev}%

\address{Department of Mathematics and Department Mechanical Engineering,
Ben-Gurion University\\ P.O.Box~653,
Beer-Sheva 84105 Israel}

\email{rsegev@bgu.ac.il}
\date{\today\\
Department of Mathematics and Department of
Mechanical Engineering,
Ben-Gurion University, Beer Sheva 84105, Israel,\\
E-mail: rsegev@bgu.ac.il, Fax:
972-8-6472814, Phone: 972-8-6477043}


\begin{abstract}
\changed{Bounds on the trace mappings
defined on the Sobolev space $\sob\ofO$ and
the space $\ldo$ of integrable stains are
obtained. Such bounds correspond to stress
concentration---the ratio between the
maximal stress in a body and the maximum of
the traction applied to its boundary.  The
analysis leading to the bounds may be
described in the mechanical context of
stress theory and stress concentration.}
%
\end{abstract}

\maketitle
\REMOVE{
\begin{flushright}
\REMOVE{
\textit{%
In memory of my friends
}\hspace*{3cm}\\{\textsc{
{Isaac Feldman, 1954--1973,\\
Amir Moses, 1954--1996,\\
Ilan Ramon (Wolferman), 1954--2003.}}}\\
\textit{%
} }
{
\aalcuin
In memory of my friends\hspace*{3.5cm}} \\
{\alcuin ISAAC FELDMAN, 1954--1973,\\
AMIR MOSES, 1954--1996,\\
ILAN RAMON (WOLFERMAN), 1954--2003. }
\end{flushright}
\vspace*{0.3cm}
}
\thispagestyle {empty}
\section{{Introduction}}
   This work
considers some mathematical aspects of
stress concentration.  Stress concentration
is a term used by engineers to indicate the
increase of stresses from some expected
values due to deviation of the geometry of
the body from an idealized simple one.  In
practice, stress concentration factors are
used by engineers to indicate the ratio
between the maximum of the stress for the
actual body under consideration and the
maximal stresses calculated for the
simplified geometry for which simple
formulae of strength-of-materials are used
traditionally.  Stress concentration factors
are usually compiled for homogeneous,
linear, isotropic, elastic solids for
various typical geometries (see for example
\cite{Peterson}). Their values are obtained
by solving the equations of elasticity
analytically, by numerical approximations
and by experimental methods.

In recent work, \cite{SCF,InvSCF}, we
formalized the notion of stress
concentration mathematically and generalized
it.  The idea in formalizing the stress
concentration factor is to regard the
stresses for the simplified geometry as
traction boundary conditions. For the
example of a bar under tension $\fc$, the
"nominal" stresses $\fc/A$, where $A$ is the
cross section area, are regarded as boundary
conditions at the ends of the bar in
agreement with the engineering notion. Thus,
for a given surface traction distribution
$\sfc$ on the boundary of the body $\gO$ and
a stress tensor $\st$ in equilibrium with
$\sfc$, the stress concentration factor is
defined by
\begin{equation}
    K_{\sfc,\st}=
 \frac{\ess\sup_{x}{\abs{\st(x)}}}
      {\ess\sup_{y}{\abs{\sfc(y)}}}\,,
      \quad
      x\in\gO,\; y\in\bndo.
\end{equation}
Here, $\abs{\st(x)}$ is a the norm of the
value of the stress at $x$.  To evaluate it
we use some norm on the vector space of
matrices.  The value of $K_{\sfc,\st}$
clearly depends on the norm used and various
norms are discussed in Section~4. Using the
essential supremum, we ignore high stresses
on sets of zero volume.

Noting that without specifying particular
constitutive relations there is a class
$\sts_\sfc$ of stress fields that are in
equilibrium with $\sfc$, we define the
optimal stress concentration factor as
\begin{equation}
    K_\sfc=\inf_{\st\in\sts_\sfc}K_{\sfc,\st}.
\end{equation}

Next, realizing that an engineer usually
does not know a-priori the nature of the
loads acting on a body exactly, the
\emph{generalized stress concentration
factor} is defined as
\begin{equation}
    K=\sup_\sfc K_\sfc
    =\sup_{\sfc}\braces{\inf_{\st\in\sts_{t}}
    \frac{\esssup_{x\in\gO}\abs{\st(x)}}%
     {\esssup_{y\in\bndo}\abs{\sfc(y)}}}
\end{equation}
where the supremum is taken over all
traction fields---essentially bounded vector
fields on $\bndo$---i.e., over all
$\sfc\in\linf\vonbdo$.  It is noted that $K$
is a purely geometric property of the body
$\gO$.

In \cite{SCF,InvSCF} we related the
generalized stress concentration factor to
the norms of the trace mappings on Sobolev
spaces and $\ld$-spaces.  It turns out that
for the formulation of equilibrium and
stress theory, particularly in the context
of stress concentration, the Sobolev space
$\sob\vono$ and the related $\ldo$ space are
especially useful.  We recall that the space
$\ldo$ contains integrable vector fields
$\vf$ such that the components of their
associated stretching, or (infinitesimal,
linear) strain,
\begin{equation}
    \eps(\vf)=\half(\nabla\vf+(\nabla\vf)^T),
    \qquad \eps(\vf)_{ij}=
 \half(\vf_{i,j}+\vf_{j,i}),
\end{equation}
are also integrable (see
\cite{TemamStrang,Temam81,Temam85,Ebobisse}).

 Assuming that $\gO$ is
an open subset of $\reals^n$ having a
$\cee2$-boundary, for both spaces one has a
well defined, linear, bounded trace mapping
$\trace$, such that for every vector field
$\vf$ defined on $\gO$, $\trace(\vf)$ is a
vector field defined on $\bndo$ satisfying
the following compatibility condition.  For
every continuous vector field $u$ defined on
the closure $\cloo$, $\trace$ acts as the
restriction to the boundary, i.e.,
\begin{equation}
    \trace(u|_\gO)=u|_\bndo.
\end{equation}

In \cite{SCF} we have shown that if we
ignore the requirement that the total force
and total torque on every subbody of $\gO$
vanish, then,
\begin{equation}
    K=\norm{\trace},\quad\text{for}\quad
     \trace\from\sob\vono\tto\lone\vonbdo.
\end{equation}
For the case where the total forces and
torques on the various subbodies do vanish,
a more detailed analysis is required (see
\cite{InvSCF}).  Letting $\rigs$ be the
finite dimensional vector space of rigid
vector fields, one has to consider the
quotient spaces $\ldo\qrigs$ and
$\lone\vonbdo\qrigs$.  For these spaces, one
can define an induced trace mapping
\begin{equation}
    \trace\qrigs\from\ldo\qrigs\tto\lone\vonbdo\qrigs,
\end{equation}
and it turns out that
\begin{equation}
    K=\norm{\trace\qrigs}.
\end{equation}

Thus, apart from the mathematical interest
in estimates on $\norm\trace$ and
$\norm{\trace\qrigs}$, such estimates are
very significant in stress analysis. It is
our objective here to estimate these
constants.  Our method of estimation has a
mechanical flavor.  In a way, it is dual to
the analysis leading to the relation between
generalized stress concentration factors and
the norms of the trace mappings.

For the particular case of the Sobolev space
$\sob\ofO$, M.~Motron obtained recently
\cite{Motron} some estimates on the bounds.
The method used here is different and is
based on the maximum principle for the
Dirichlet problem.  The bounds we obtain
give a concrete estimate for $\norm\trace$.
Subsequently, we extend our method to the
space of $\ld$-fields.
%
Specifically,
recalling that the $\ld$-norm is given by
\begin{equation}
    \norm\vf_\ld=\sum_i\norm{\vf_i}_\lone+
            \sum_{i,m}\norm{\eps(\vf)_{im}}_\lone,
\end{equation}
we obtain bounds on the constants $A$ and
$B$ such that
\begin{equation}
\ino{\bdry\gO}\abs{\trace(\vf)}\les
   A\ino\gO\abs\vf\dV+B\ino\gO\abs{\eps(\vf)}
\end{equation}
for all $\ld$-fields $\vf$ (so that
$\max\{A,B\}$ bounds $\norm\trace$).

Section 2 presents the results and methods
of \cite{InvSCF} relating stress
concentration to the norm of the trace
mapping.  The last subsection discusses the
simplification to the case of sourceless
vector fields (rather then stress tensors)
which may serve as motivation for studying
the norm of the trace mapping for the
Sobolev space $\sob(\gO)$.  In Section 3 we
introduce the basic method for obtaining the
bounds on the trace mapping for the Sobolev
space using harmonic vector fields. The main
result of this section is Theorem
\ref{thm:ronen}.  As background material for
the discussion of the bounds for the trace
mapping on $\ldo$, Section 4 presents some
standard results on norms of matrices. These
are significant in the mechanical context as
they are used on the space of stress
matrices.  For example, a yield criterion is
usually a seminorm on the space of stress
matrices.  Section 5 presents additional
preparatory material---the optimal boundary
values for the stresses for a given boundary
traction field. \changed{Section 6 studies
the bounds on the norm of the trace mapping
for $\ldo$ using the method of harmonic
tensor fields and the central result is
Theorem~\ref{thm:BestConstLD}. Finally, the
concluding remarks of Section 7 discuss the
mechanical interpretation of the preceding
analysis.}


\section{Generalized Stress Concentration and\\
    the Norm of the Trace Mapping}

\subsection{Basic definitions and notation}
We consider an open set
$\gO\subset\reals^3$, where in some
sections, the presentation is in the setting
of $\reals^n$. We assume that $\gO$ is
bounded and that it has a $\cee2$-boundary.
(The results hold for less restrictive
assumptions.)  We will use the index
summation convention for repeated indices
and subscripted comma followed by an index
will indicate partial differentiation with
respect to the corresponding variable.

A vector field on $\gO$ is interpreted
physically as a virtual velocity field on
the body or alternatively as a field of
virtual infinitesimal displacements. A
\emph{rigid} field in $\rthree$ is a vector
field of the form
\[
 \vf(x)=a+b\times x,\quad
a,b\in\rthree.
\]
 Clearly, rigid fields may
be restricted to subsets of $\rthree$. We
denote the space of rigid fields by $\rigs$
and it is a 6-dimensional vector space.

The following definitions and results
concerning $\ld$-fields are due to Temam and
Strang \cite{TemamStrang,Temam81,Temam85}.
Given an integrable vector field $\vf$ on
$\gO$, we consider the corresponding
\emph{stretching} (\emph{linear strain})
field $\eps(\vf)$ defined by\begin{equation}
\eps(\vf)_{im}=\half(\vf_{i,m}+\vf_{m,i}),\end{equation}
where comma implies distributional
derivative relative to the corresponding
spatial coordinate. The integrable vector
field $\vf$ is of \emph{integrable
stretching}, or $\vf\in\ldo$, if the
components of the corresponding stretching
are also integrable over $\gO$. On the
vector space $\ldo$ of integrable
stretchings it is natural to use the norm
\begin{equation}
\norm\vf=\norm{\vf}_{\ld}=
\norm{\vf}_1+\norm{\eps(\vf)}_1,\label{eq:ld-norm-def}
\end{equation}
where $\norm\cdot_p$ indicates the
$L^p$-norm.  With this norm, $\ldo$ is a
Banach space and we have a continuous and
linear\begin{equation}
\eps\from\ldo\tto\lone\strono.\end{equation}
A basic theorem whose classical version is
due to Liouville (see \cite[ pp.
18-19]{Temam85}) states:

\begin{prop}
$\kernel(\eps)=\rigs$.
\end{prop}
Let $\vs$ be a Banach space of velocity
fields. In the discussion below $\vs$ will
be either $\ldo$ or $\lone\vonbdo$. We refer
to an element $\dis\in\vs\qrigs$ as a
\emph{distortion}. We have the natural
projection mapping onto the quotient
space
\begin{equation}
\pq\from\vs\tto\vs\qrigs\end{equation}
 and the induced norm in $\vs\qrigs$ is given by\begin{equation}
\norm{\dis}=\inf_{\vf\in\dis}\norm{\vf},\quad\text{or}\quad\norm{[\vf]}=\inf_{\rig\in\rigs}\norm{\vf+\rig}.\end{equation}

\begin{prop}
\label{pro:ProjRigs}For both
$\vs=\lone\vonbdo$ and $\vs=\ldo$ there are
continuous and linear projection mappings
\begin{equation}
\pr\from\vs\tto\rigs.\end{equation} For
$\vf\in\vs$, $\pr(\vf)=a+b\times x$ is given
by\begin{equation}
a=\frac{1}{\abs{U}}\ino{U}\vf,\quad
b=I^{-1}\paren{\ino{U}x\times\vf},\label{eq:proj-rig}\end{equation}
where $U=\bndo$ for $\vs=\lone\vonbdo$,
$U=\gO$ for $\vs=\ldo$, $\abs{U}$ is the
Hausdorff measure of $U$, and
$I_{im}=\int_{U}(x_{k}x_{k}\delta_{im}-x_{i}x_{m})$
is the moment of inertia of $U$.
\end{prop}
The space $\ldo$ has the following
properties (see \cite{Temam85}).

\begin{description}
\item [Approximation]The restrictions of fields in $\cee{\infty}\vonclo$
to $\gO$ are dense in $\ldo$.
\item [Extensions]There is a continuous linear extension operator $E\from\ldo\to\ld(\rthree)$.
\item [Regularity]If $\vf$ is any distribution on $\gO$ whose corresponding
stretching is $\lone$, then
$\vf\in\lone\vono$.
\item [Trace mapping]There is a unique linear,
 surjective, continuous trace mapping
 \begin{equation}
\trace\from\ldo\tto\lone\vonbdo
\end{equation}
such that $\trace(\vf|_\gO)=\vf|_\bndo$ for
all continuous vector fields $\vf$ defined
on $\gO$.
\item [Distortions]On $\ldo\qrigs$, \begin{equation}
\norm{\dis}_{\eps}=\norm{\eps(\dis)}_1\label{eq:EqNormDistortions}\end{equation}
 is a norm that is equivalent to the quotient norm described above.
Thus, there is a constant $C\ofO$ (depending
on $\gO$ only) such that for every
$\vf\in\ldo$
\begin{equation}
\inf_{r\in\rigs}\norm{\vf+\rig}\les
C(\gO)\norm{\eps(\vf)}_1.
\end{equation}
The infimum is attainable, i.e., for each
$\vf\in\ldo$, there is a rigid motion
$\rig_0$ satisfying
\begin{equation}
    \norm{\vf+\rig_0}=\inf_{r\in\rigs}\norm{\vf+\rig}
     \les
      C(\gO)\norm{\eps(\vf)}_1
\end{equation}

\item [Equivalent~norm]If $p$ is a seminorm on
$\ldo$ such that $p(\rig)=0$ implies that
$\rig=0$ for every $\rig\in\rigs$ (so $p$ is
a norm on $\rigs$), then,\begin{equation}
p(\vf)+\norm{\eps(\vf)}_1\end{equation} is a
norm on $\ldo$ that is equivalent to the
original norm. In particular, using the
trace mapping one can take\begin{equation}
p(\vf)=\norm{\trace(\vf)}_\onebdo,\end{equation}
so the following is a norm that is
equivalent to original
(\ref{eq:ld-norm-def})\begin{equation}
\norm{\vf}_{\times}=\norm{\trace(\vf)}_\onebdo+\norm{\eps(\vf)}_1.\end{equation}
Furthermore, one may use the projection
$\pr\from\lone\vonbdo\to\rigs$ as in
Proposition~\ref{pro:ProjRigs} and a norm
$\norm{\cdot}_{\rigs}$ to obtain the
equivalent norm\begin{equation}
\norm{\vf}_{\oplus}=\norm{\pr\comp\trace(\vf)}_{\rigs}+\norm{\eps(\vf)}_1.\label{eq:EqLDNorm}\end{equation}

\end{description}

Forces are regarded as elements of the dual
spaces to the corresponding spaces of
virtual velocities.  So for a generic space
of velocities $\vs$, a force $\fc$ will be a
member of $\vs^*$.  The evaluation
$\fc(\vf)$ is interpreted as virtual work,
or virtual power, performed by the
generalized force for the corresponding
generalized velocity.

In case the space $\vs$ of velocities is an
$L^p$-space, $1\les{p}\les\infty$, (e.g.,
$\lone\vonbdo$) a force may be represented
by an element of the corresponding $L^q$
space with $q=p/(p-1)$ through integration.
We will use the same symbol for the force
and its representing field.  For example,
for $\sfc\in\lone\vonbdo^*=\linf\vonbdo$ we
have
\begin{equation}
    \sfc(\vf)=\ino\bndo\sfc\cdot\vf.
\end{equation}

A force $\fc\in\vs^{*}$ acting on a body is
\emph{equilibrated} if $\fc(\rig)=0$ for all
$\rig\in\rigs$. An equilibrated force $\fc$
is of the form $\fc=\pq^{*}(\fc_{0})$ for
some $\fc_{0}\in(\vs\qrigs)^{*}$.
Furthermore, $\pq^{*}$ is norm preserving,
i.e.,\begin{equation}
\norm{\pq^*(\fc_{0})}=\norm{\fc_{0}},\label{eq:PiDualNormPreserving}\end{equation}
 so one can usually identify an equilibrated $\fc$ with $\fc_{0}$.

\subsection{Generalized stress concentration
factors and norms of trace mappings.}

The central mathematical object that we find
suitable for formulating the continuum
mechanics problem, particularly, those
notions related to stress concentration is
$(\ldo\qrigs)^{*}$---the dual to the space
of $\ld$-distortions. Specifically, as
described below and in further detail in
\cite{InvSCF}, elements of this space may be
represented by essentially bounded stress
fields and on the other hand, the dual of
the trace mapping associates an element of
$(\ldo\qrigs)^{*}$ with any equilibrated
boundary traction field.

Consider the composite mapping
$\pq\comp\trace\from\ldo\to\lone\vonbdo\qrigs$.
It is noted that for any $\rig\in\rigs$,
$\trace(\rig)$ is a rigid motion on $\bndo$,
hence,
$\pq\comp\trace(\vf+\rig)=\pq\comp\trace(\vf)+\pq\comp\trace(r)=\pq\comp\trace(\vf)$.
Thus, we have a well defined
mapping\begin{equation}
\trace\qrigs\from\ldo\qrigs\tto\lone\vonbdo\qrigs,\end{equation}
given by
$\trace\qrigs(\dis)=\pq\comp\trace(\vf)$,
for some $\vf\in\dis$.
Clearly,
\begin{equation}
\pq\comp\trace=(\trace\qrigs)\comp\pq.
\end{equation}
In addition,
\begin{equation}
\norm{\trace(\vf)+\rig}_1
=\norm{\trace(\vf+\rig)}_1\les\norm{\trace}\norm{\vf+\rig}.
\end{equation}
Hence,
\begin{equation}
\norm{[\trace(\vf)]}
=\inf_{\rig\in\rigs}\norm{\trace(\vf)+\rig}_1\les\norm{\trace}\inf_{r\in\rigs}\norm{\vf+\rig}
=\norm{\trace}\norm{[\vf]}
\end{equation}
and we conclude that $\trace\qrigs$ is
indeed bounded and
\begin{equation}\label{eq:TraceQrigsBounded}
\norm{\trace\qrigs}\les\norm{\trace}.
\end{equation}

The dual mapping
$\trace^{*}\from\linf\vonbdo\to(\ldo)^{*}$
may now be applied to traction fields and
$(\trace\qrigs)^{*}\from(\lone\vonbdo\qrigs)^{*}\to(\ldo\qrigs)^{*}$
may be applied to equilibrated boundary
traction fields to give $\ld$-forces and
equilibrated $\ld$-forces respectively.
Clearly, \begin{equation}
\pq^{*}\comp(\trace\qrigs)^{*}=\trace^{*}\comp\pq^{*}.\label{eq:CommutForces}\end{equation}

\begin{prop}
The mappings $\trace^{*}$ and
$(\trace\qrigs)^{*}$ are injective.
\end{prop}
\begin{proof}
The mapping $\trace^{*}$ is injective
because $\trace$ is continuous and
surjective. As the quotient space projection
$\pq\from\lone\vonbdo\to\lone\vonbdo\qrigs$
is also continuous and surjective, the same
argument applies to $\trace\qrigs$.
\end{proof}
\begin{remk}
\label{rem:ChoicOfNorms}
Henceforth, we will use the equivalent norm
$\norm{\cdot}_{\eps}$ as in
(\ref{eq:EqNormDistortions}) on
$\ldo\qrigs.$ We will also use the
equivalent norm $\norm{\cdot}_{\oplus}$ on
$\ldo$ as in (\ref{eq:EqLDNorm}). This
implies that the quotient norm on
$\ldo\qrigs$ is actually equal (not only
equivalent) to the norm induced by the
strain. That is,\begin{equation}
\norm{[\vf]}=\norm{[\vf]}_{\eps}=\norm{\eps(\vf)}_1.\end{equation}
It is noted that the last equality, to be
used below frequently, is independent of the
choice of a particular projection
$\pr\from\ldo\to\rigs$.
\end{remk}
\begin{prop}
\label{pro:RepByStresses}
Any $\ebdfc\in(\ld\ofO\qrigs)^{*}$ is
represented by some symmetric stress field
$\st\in\linf\strono^{*}$, in the form
\begin{equation}
\ebdfc=\eps^{*}(\st),
\end{equation}
where
$\eps^{*}\from\linf\strono\to(\ldo\qrigs)^{*}$
is the dual mapping to
\begin{equation}
\eps\from\ldo\qrigs\tto\lone\strono.
\end{equation}
In addition, for the dual norm
$\norm{\cdot}^{\eps}$ on $(\ldo\qrigs)^{*}$,
we have
\begin{equation}
\norm{\ebdfc}^{\eps}=\inf_{\st,\,\ebdfc
=\eps^{*}(\st)}\norm{\st}_{\infty},
\end{equation}
and the infimum is attained for some
$\hat{\st}\in\linf\strono$,
i.e.,
\begin{equation}
\norm{\ebdfc}^{\eps}=\norm{\hat{\st}}_\infty.
\end{equation}

\end{prop}

\begin{proof}
Using the duality
$\lone\strono^*=\linf\strono$, the assertion
follows from the fact that
$\eps\from\ldo\qrigs\to\lone\strono$ is a
norm-preserving (by our choice of norm as in
Remark~\ref{rem:ChoicOfNorms}), linear
injection and using the Hahn-Banach theorem
(see \cite{InvSCF} for the details).
\end{proof}
\begin{corl}
\label{cor:EquilibriumCondition}
Let $\sfc\in\linf\vonbdo$ be any
equilibrated traction field so there is a
$\sfc_{0}\in(\lone\vonbdo\qrigs)^{*}$ such
that $\sfc=\pq^{*}(\sfc_{0})$. Then, there
exists some stress field
$\st\in\linf\strono$ such that
\begin{equation}\label{eq:EquilCond}
(\trace\qrigs)^{*}(\sfc_{0})
  =\eps^{*}(\st),\quad\text{and}\quad\trace^{*}(\sfc)
  =\pq^{*}\comp\eps^{*}(\st).
\end{equation}
In addition,
\begin{equation}
\norm{(\trace\qrigs)^{*}(\sfc_{0})}
=\inf_{\eps^{*}(\st)
=(\trace\qrigs)^{*}(\sfc_{0})}
 \,\,\norm{\st}_{\infty},
\end{equation}
and
\begin{equation}
\norm{\trace^{*}(\sfc)}
=\inf_{\pq^{*}\comp\eps^{*}(\st)
 =\trace^{*}(\sfc)}\,\,\norm{\st}_{\infty}.
\end{equation}

\end{corl}
\begin{proof}
We make repetitive use of
(\ref{eq:CommutForces}),
(\ref{eq:PiDualNormPreserving}) and
Proposition (\ref{pro:RepByStresses}). For
example,
\begin{equation}\begin{split}
\norm{\trace^{*}(\sfc)} & =\norm{\trace^{*}\comp\pq^{*}(\sfc_{0})}\\
 & =\norm{\pq^{*}\comp(\trace\qrigs)^{*}(\sfc_{0})}\\
 & =\norm{(\trace\qrigs)^{*}(\sfc_{0})}\\
 & =\inf_{\pq^{*}\comp\eps^{*}(\st)=\pq^{*}\comp(\trace\qrigs)^{*}(\sfc_{0})}\,\,\norm{\st}_{\infty}\\
 & =\inf_{\pq^{*}\comp\eps^{*}(\st)=\trace^{*}\comp\pq^{*}(\sfc_{0})}\,\,\norm{\st}_{\infty}\\
 & =\inf_{\pq^{*}\comp\eps^{*}(\st)=\trace^{*}(\sfc)}
 \,\,\norm{\st}_{\infty}.
 \end{split}
 \end{equation}
\end{proof}
\begin{remk}
\label{rem:EquilCond}
The conditions (\ref{eq:EquilCond}) are
equivalent to the principle of virtual
work---a weak form of the equations of
equilibrium---of continuum mechanics (as it
is assumed throughout that the body forces
vanish). For example,
$\trace^{*}(\sfc)=\pq^{*}\comp\eps^{*}(\st)$,
implies
\begin{align}
 \trace^{*}(\sfc)(\vf) &
=(\pq^{*}\comp\eps^{*})(\st)(\vf),
\end{align}
so that for any
$\vf\in\ldo$,\begin{equation}
\sfc(\trace(\vf))=\st(\eps\comp\pq(\vf))=\st(\eps(\vf)).\end{equation}
Hence, for a vector field $\vf$ that is the
restriction of a differentiable field on
$\cloo$,\begin{equation}
\ino{\bndo}\sfc_{i}\vf_{i}=\ino{\gO}\st_{ij}\eps(\vf)_{ij}.\end{equation}

\end{remk}
\begin{thm}
\label{thm:KEqualNormTrace}
Let,
\begin{equation}
\norm{\trace\qrigs}=\sup_{\dis\in\ldo\qrigs}\frac{\norm{(\trace\qrigs)(\dis)}}{\norm{\dis}},\end{equation}
be the norm of the trace mapping for
distortions. Then, the stress concentration
factor $K$ for the boundary traction
problem, satisfies\begin{equation}
K=\norm{\trace\qrigs}.\end{equation}
Specifically,
\begin{equation}
K=\sup_{\vf\in\cee{\infty}\ofclo}
 \frac{\inf_{\rig\in\rigs}\norm{\vf|_\bndo+r}_1}
  {\norm{\eps(\vf)}_1}\label{eq:ApprKbySmooth}
  \end{equation}

\end{thm}
\begin{proof}
We have the standard \begin{equation}
\norm{\trace\qrigs}=\norm{(\trace\qrigs)^{*}}.\end{equation}
However,\begin{align}
\norm{(\trace\qrigs)^{*}} & =\sup_{\sfc_{0}}\frac{\norm{(\trace\qrigs)^{*}(\sfc_{0})}^{\eps}}{\norm{\sfc_{0}}}\\
 & =\sup_{\sfc_{0}\in(\vonbdo\qrigs)^{*}}\braces{\frac{1}{\norm{\sfc_{0}}}\,\,\inf_{\eps^{*}(\st)=(\trace\qrigs)^{*}(\sfc_{0})}\norm{\st}_{\infty}}\\
 & =\sup_{\sfc\in\image\pq^{*}}\braces{\frac{1}{\norm{\sfc}_{\infty}}\,\,\inf_{\pq^{*}\comp\eps^{*}(\st)=\trace^{*}(\sfc)}\,\,\norm{\st}_{\infty}},\end{align}
where we used Proposition
(\ref{pro:RepByStresses}) and Corollary
(\ref{cor:EquilibriumCondition}). The
condition
$\pq^{*}\comp\eps^{*}(\st)=\trace^{*}(\sfc)$
is equivalent to the condition that the
stress field $\st$ is in equilibrium with
$\sfc$, i.e., $\st\in\sts_{\sfc}$, by Remark
(\ref{rem:EquilCond}). The expression
(\ref{eq:ApprKbySmooth}) simply uses the
definition of $\norm{\trace\qrigs}$ and the
fact that $\cee{\infty}\ofclo$ is dense in
$\ldo$.
\end{proof}

\subsection{The scalar case and the trace mapping on the Sobolev space $\sob\ofO$}
The previous discussion is simplified
considerably if we consider scalar fields
$\fun\in\sob\ofO$ instead of the vector
fields $\vf\in\ldo$. In this case the
boundary data is also a scalar field, the
analog of a stress is a vector field, and
the vector space $\rigs$ of rigid motions is
replaced by the real numbers.

Consider a sourceless vector field $\st$ on
$\gO$ that satisfies boundary conditions for
its boundary flux, i.e.,
\begin{alignat}{2}
\st_{i,i} & =0,\quad & \text{in }\gO,\label{eq:field}\\
\st_{i}\nor_{i} & =\sfc,\quad & \text{on
}\bndo,\label{eq:boundaryCond}
\end{alignat}
for some given essentially bounded
$\sfc\from\bndo\to\reals$. Physically, $\st$
may be thought of as a material flow field
(say for an incompressible flow) for a given
flux density $\sfc$ on the boundary.
Alternatively, $\st$ may be thought of as a
heat flow field where there are no heat
sources in $\gO$ so $\sfc$ is the given heat
flux on the boundary; or $\st$ may be the
electric displacement field and $\sfc$ is
the charge density on the boundary.

The weak formulation of the problem is
\begin{equation}
\ino{\gO}\st_{i}\fun_{,i}
 =\ino{\bndo}\sfc\fun.\label{eq:weakForm}
 \end{equation}
(The test function $\fun$ may be thought of
as a potential field in the electrostatic
example or as the reciprocal of the
temperature in the heat transfer example.)

The analog of the stress concentration
factor is then
\begin{equation}
K=\sup_{\sfc\in\linf\ofbdo}K_{\sfc}
 =\sup_{\sfc\in\linf\ofbdo}
  \inf_{\st\in\sts_{t}}\frac{\norm{\st}_{\infclo}}
  {\norm{\sfc}_{\infbdo}}
   =\sup_{\sfc\in\linf\ofbdo}\inf_{\st\in\sts_{t}}
    \frac{\esssup_{x\in\gO}\abs{\st(x)}}%
     {\esssup_{y\in\bndo}\abs{\sfc(y)}}.
\end{equation}
We will refer to $K$ as the
\emph{generalized field concentration
factor}.  Thus for example, for the
interpretation of $\st$ as a flow field, for
a given flux density $\sfc$, $K_t$ will be
the smallest ratio between the maximal
magnitude of the velocity and the maximal
value of the given boundary flux.

Thus, $\fun$ may be regarded as an element
of the Sobolev space $\sob\ofO$ so
\begin{equation}
\norm{\fun}_{\sob}=\norm{\fun}_1+\norm{\nabla\fun}_1.\end{equation}
 On the Sobolev space, the trace mapping\begin{equation}
\trace\from\sob\ofO\to\lone\ofbdo\end{equation}
is well defined as expected (see
\cite{Adams}). Our assumption that there are
no sources in $\gO$ implies that
$\int_{\bndo}t=0$ which is equivalent to
considering $\sob\ofO\qr$ and $\trace\qr$.
Thus, the analog of Theorem
(\ref{thm:KEqualNormTrace}) will be for the
scalar case \begin{equation}
K=\norm{\trace\qr},\quad\trace\from\sob\ofO\tto\lone\ofbdo.\end{equation}


\section{Bounds on the $\sob$-Trace Operator}
 \label{sobTrace}
\subsection{The bounds obtained using normal vector fields}
In this section we consider bounds for the
trace operator
\begin{equation}
\trace\from\sob\ofO\tto\lone(\bndo).
\end{equation}
In particular, we are looking for bounds $A$
and $B$ satisfying
\begin{equation}
\ino\bndo\abs\fun\les
 A\ino\gO\abs{\nabla\fun}
  +B\ino\gO\abs\fun
\end{equation}
for every $\fun\in\sob\ofO$.

Let $\nvf$ be any $\cee1$-vector field on
$\cloo$ and $\psi$ the restriction to
$\cloo$ of a $\sob$-function defined in an
open neighborhood of $\cloo$.  Then,
\begin{equation}
\ino\gO\nvf_i\psi_{,i}=\ino\gO\paren{\nvf_i\psi}_{,i}
          -\ino\gO\nvf_{i,i}\psi
\end{equation}
implies using the Gauss-Green theorem that
\begin{equation}\label{eq:pvw}
\ino\bndo\nvf_i\nor_i\psi=\ino\gO\nvf_i\psi_{,i}
           +\ino\gO\nvf_{i,i}\psi,
\end{equation}
where $\nor$ is the outwards pointing unit
normal to $\bndo$.
\begin{defn}\label{def:NormalVf}
The vector field $\nvf$ on $\cloo$ will be
referred to as a normal field if the
following conditions hold.
\begin{enumerate}
\item[(\textit{i})]  $\nvf(y)=\nor(y)$ for all
$y\in\bndo$.
\item[(\textit{ii})]  $\abs\nvf(x)\les1$ for
all $x\in\gO$, where here and in the rest of
this section we use the Euclidean norm for
elements of $\reals^n$ so
$\abs\nvf=\sqrt{\nvf_i\nvf_i}$.
\end{enumerate}
\end{defn}
The existence of normal vector fields is
discussed in some detail below.  We will use
$\nvfs\ofO$ for the collection of all normal
vector fields. For a normal field,
Equation~(\ref{eq:pvw}) assumes the form
\begin{equation}\label{eq:pvwNorVF}
   \ino\bndo\psi=\ino\gO\nvf_i\psi_{,i}
           +\ino\gO\nvf_{i,i}\psi.
\end{equation}
Given a $\cee1$ mapping $\fun$ on $\cloo$,
the distributional derivatives
$\abs{\fun}_{,i}$ of its absolute value
$\abs{\fun}(x)=\abs{\fun(x)}$ are clearly
integrable and hence, $\abs\fun$ is $\sob$.
Rewriting Equation~(\ref{eq:pvwNorVF}) for
$\psi=\abs\fun$ we obtain
\begin{equation}
    \ino\bndo\abs\fun=\ino\gO\nvf_i\abs{\fun}_{,i}
           +\ino\gO\nvf_{i,i}\abs\fun.
\end{equation}
We now estimate each of the integrals on the
right hand side.
\begin{equation}
\begin{split}
    \ino\gO\nvf_i\abs{\fun}_{,i}&
       \les\ino\gO\abs\nvf\abs{\abs{\fun}_{,i}}\\
    &\les\ino\gO\abs{\nabla\fun},
\end{split}
\end{equation}
where we used Definition
(\ref{def:NormalVf}.\emph{ii}) and
$\abs{\abs{\fun}_{,i}}=\abs{\fun_{,i}}\les\abs{\nabla\fun}$.
Also
\begin{equation}
        \ino\gO\nvf_{i,i}\abs\fun
            \les\max_{\cloo}\braces{\abs{\nvf_{i,i}}}\ino\gO\abs\fun.
\end{equation}
It follows that
\begin{equation}
    \ino\bndo\abs\fun\les
     \ino\gO\abs{\nabla\fun}+
      \max_{x\in\cloo}\braces{\abs{\nvf_{i,i}(x)}}\ino\gO\abs\fun.
\end{equation}
{This inequality is clearly exact} and
recalling the definition of normal vector
fields we have
\begin{thm}
Let $\gO$ be a bounded open set in
$\reals^n$ having a $\cee2$-boundary and set
\begin{equation}
    B(\gO)=\inf_{\nvf\in\nvfs\ofO}
     \braces{\max_{x\in\cloo}\braces{\abs{\nvf_{i,i}(x)}}}
     =\inf_{\nvf\in\nvfs\ofO}\braces{\norm{\nvf_{i,i}}_\infclo}.
\end{equation}
Then, the following exact inequality holds
\begin{equation}
    \norm\fun_\onebdo=\ino\bndo\abs\fun\les
     \ino\gO\abs{\nabla\fun}+B(\gO)\ino\gO\abs\fun
      =\norm{\nabla\fun}_1+B(\gO)\norm\fun_1.
\end{equation}
\end{thm}
\subsection{Estimation using harmonic normal fields}
\label{subs:EstHarmVF}
We now consider the existence of normal
vector fields.  By the assumption that
$\bndo$ is $\cee2$, $\nor$ is a
$\cee1$-vector field on $\bndo$ and  we can
extend it to a vector field $\nvf$ with
$\abs{\nvf(x)}\les1$ for all $x\in\cloo$
(see for example Theorem 3.6.2 in
\cite{Ziemer}). Furthermore, we can require
that the extension is harmonic in the
following sense. For a vector field $\nvf$,
we use $\Delta\nvf$ for the vector field
$(\Delta\nvf)_j=(\nvf_j)_{,ii}$.  The field
$\nvf$ is \emph{harmonic} in $\gO$ if
\begin{equation}
    \Delta\nvf=0,\quad\text{for all\
    }x\in\gO.
\end{equation}

\begin{thm}\label{thm:harVF}
Let $\gO$ be a bounded open subset of
$\reals^n$ such that $\bndo$ is $\cee2$.
Then, there exists a unique normal vector
field $\nvf_0\in\nvfs\ofO$ which is
harmonic. In addition, for the supremum of
the divergence,
$\nabla\cdot\nvf_0=\nvf_{0i,i}$, we have
\begin{equation}\label{eq:maxPrDiv}
    \norm{\nabla\cdot\nvf_0}_\infclo
     =\norm{\nabla\cdot\nvf_0}_\infbdo.
\end{equation}
\end{thm}
\begin{proof}
For any fixed $j$, we have a classical
Dirichlet problem
\begin{equation}
    \Delta\nvf_j=0\quad
    \text{in}\quad\gO,\qquad
    \nvf_j=\nor_j \quad \text{on}\quad\bndo.
\end{equation}
Given our smoothness assumption on $\bndo$,
there is a unique solution $\nvf_{0j}$ to
each such boundary value problem and we
obtain the harmonic vector field $\nvf_0$.

For $\abs{\nvf_0}^2=\nvf_{0j}\nvf_{0j}$ we
have,
\begin{equation}
\begin{split}
    \Delta(\nvf_{0j}\nvf_{0j})&=(\nvf_{0j}\nvf_{0j})_{,ii}\\
        &=(2\nvf_{0j,i}\nvf_{0j})_{,i}\\
        &=2\nvf_{0j,ii}\nvf_{0j}+2\nvf_{0j,i}\nvf_{0j,i}.
\end{split}
\end{equation}
In the last line, the first term vanishes
because $\nvf_{0j}$ is harmonic and hence
$\Delta(\abs{\nvf_0}^2)\ges0$.  We conclude
that $\abs{\nvf_0}^2$ is subharmonic in
$\gO$.  By the maximum principle for
subharmonic functions
\begin{equation}
    \max_{x\in\cloo}\abs{\nvf_0(x)}^2=
    \max_{y\in\bndo}\abs{\nvf_0(y)}^2=
    \max_{y\in\bndo}\abs{\nor(y)}^2=1.
\end{equation}
Thus, in addition to the boundary
conditions, $\nvf_0$ satisfies the condition
(\ref{def:NormalVf}.\textit{ii}), and so
$\nvf_0\in\nvfs\ofO$.

Next, we note that for the harmonic
$\nvf_0$,
\begin{equation}
\begin{split}
    \Delta(\nabla\cdot\nvf_0)&=(\nvf_{0j,j})_{,ii}\\
       &=(\nvf_{0j,ii})_{,j}\\
       &=0.
\end{split}
\end{equation}
Thus, $\nabla\cdot\nvf_0$ is also harmonic
in $\gO$ so by the maximum principle
Equation (\ref{eq:maxPrDiv}) holds.
\end{proof}
It turns out that the harmonic vector field
of Theorem (\ref{thm:harVF}) plays an
important role in the computation of
$B\ofO$, the second constant of the sobolev
$\sob\ofO$ trace inequality.  Continuing to
use $\nvf_0$ for the unique harmonic normal
vector field we have
\begin{thm}\label{thm:ronen}
The Sobolev constant $B\ofO$ is given by
\begin{equation}
    B\ofO=
     \inf_{\nvf\in\nvfs\ofO}\norm{\nabla\cdot\nvf}_\infclo
     =\norm{\nabla\cdot\nvf_0}_{\infbdo}.
\end{equation}
\end{thm}
\begin{proof}
Let $(\nvf_m)$, $m\in\mathbb{N}$, be a
sequence of normal vector fields such that
\begin{equation}
    \lim_{m\to\infty}\norm{\nabla\cdot\nvf_m}=B\ofO.
\end{equation}
Using normal tubular neighborhoods (e.g.,
\cite[p.~110]{Hirsch}), there is a
$\delta>0$ such that we can parameterize an
open neighborhood $V$ of $\bndo$ in $\gO$ by
\begin{equation}
    (y,z)\in\bndo\times[0,\delta),
\end{equation}
where for each $x\in V$, $y(x)$ is a unique
point on the boundary such that $x$ is on
the line through $y$ which is normal to the
boundary, and $z$ is the distance to the
boundary along the normal line where $x$ is
situated. For each $m$, let $V_m$ be the
open set
\begin{equation}
    V_m=\braces{x\in\cloo:
     {y(x)\in\bndo,z(x)<\delta/m}},
\end{equation}
and let $\gO_m=\cloo-\cl{V}_m$ so
\begin{equation}
    \bndo_m=\braces{x\in\cloo:
     {y(x)\in\bndo,z(x)=\delta/m}}.
\end{equation}

Now for each $m$ we construct the harmonic
lifting (cf. \cite[p.~24]{Gilbarg})
$\cl\nvf_m$ of $\nvf_m$ as follows.  Let
$\nvf_{0m}$ be the solution of the Dirichlet
problem on $\gO_m$ with the boundary
conditions $\nvf_{0m}(x)=\nvf_m(x)$ for
$x\in\bndo_m$.  Set
\begin{equation}
    \cl\nvf_m(x)=\begin{cases}
    \nvf_{0m}(x),\quad&\text{for } x\in\gO_m,\\
    \nvf_{m}(x),\quad&\text{for } x\in{\cl V}_m.
    \end{cases}
\end{equation}
By the maximum principle,
\begin{equation}
    \norm{\cl\nvf_{mi}}_\infclo
          \les\norm{\nvf_{mi}}_\infclo
  \quad\text{and}\quad
     \norm{\nabla\cdot\cl\nvf_{m}}_\infclo
           \les\norm{\nabla\cdot\nvf_{m}}_\infclo,
\end{equation}
so
\begin{equation}
    \lim_{m\to\infty}\norm{\nabla\cdot\cl\nvf_m}_\infclo
           =B\ofO.
\end{equation}
In addition, as the various
$\norm{\nvf_{mi}}_\infclo$ are bounded by 1,
the same applies to $\cl\nvf_{mi}$ and the
sequence $\cl\nvf_{mi}$ is uniformly
bounded. Thus, (using a standard normal
family argument) by Ascoli's theorem, it has
a subsequence that converges uniformly to a
limit continuous normal field $\cl\nvf$. On
any compact subset of $\gO$ this gives a
uniformly convergent sequence of harmonic
functions whose limit is then also harmonic.
Thus, $\cl\nvf$ is harmonic. Also, the limit
$\cl\nvf$ satisfies the conditions of
Definition (\ref{def:NormalVf}) and so it is
a normal vector field.  Finally, by the
uniqueness of the solution to the Dirichlet
problem $\cl\nvf=\nvf_0$.
\end{proof}

\begin{remk}
Within the framework of the $AB$-program in
geometric analysis, \cite{Druet}, M.~Motron
proves in \cite{Motron}, using different
methods, the following two theorems for
bounds on the trace mapping on $\sob\ofO$.
\begin{enumerate}
    \item For any $\eps>0$, there exists a
    $B_\eps$ such that for any
    $\fun\in\sob\ofO$,
    \begin{equation}\label{eq:MotronFirstA}
        \ino\bndo\abs\fun
        \les(1+\eps)\ino\gO\abs{\nabla\fun}
           +B_\eps\ino\gO\abs\fun.
\end{equation}
    \item Assuming that $\gO$ is a connected
    bounded open subset of $\reals^n$ whose
    boundary is piecewise $\cee1$, there
    exists $A>0$ such that for any
    $\fun\in\sob\ofO$,
    \begin{equation}
        \ino\bndo\abs\fun\les
        A\ino\gO\abs{\nabla\fun}+
         \frac{\abs\bndo}{\abs\gO}\ino\gO\abs\fun.
\end{equation}
\end{enumerate}
We note that even if we had the solution to
the $AB$-program, we would not have a
concrete bound on $\norm\trace$. What we
sought were simultaneous bounds on both $A$
and $B$.

\changed{
In addition, for a normal vector field
\begin{equation}
    \sup_{x\in\cloo}\abs{\nabla\cdot\nvf}\abs\gO
    \ges\ino\cloo\nabla\cdot\nvf
    =\ino\bndo\nvf\cdot\nor
    =\abs\bndo,
\end{equation}
so
\begin{equation}
    \inf_{\nvf\in\nvfs}\norm{\nabla\cdot\nvf}_\infclo
     \ges\frac{\abs\bndo}{\abs\gO}\,.
\end{equation}
However, the inequality is not exact and
equality is not attainable even for the
harmonic normal vector field.  (Think of two
circles connected by a narrow neck of width
$t$. Across the narrow neck,
$\nabla\cdot\nvf$ has to be of order $1/t$
in order to satisfy the boundary conditions.
The values of $\abs\gO$ and $\abs\bndo$ are
not significantly different from those of
the two circles.)  \changedd{Thus, Motron's
bound $B$ is smaller than the value obtained
here. On the other hand, the bounds $A$ and
$B$ we obtain are exact in the sense that
you cannot lower $B$ without increasing $A$.
}}
\end{remk}

\section{Norms of Symmetric Tensors}
 \label{NormsTensors}
As noted earlier, the value of the stress
concentration factor depends on the norm we
choose to use on the space of stress
matrices.  Thus, for the sake of
completeness, we review below some
elementary properties of norms of symmetric
matrices on $\reals^n$.  We will denote the
norm of a matrix $\lin$ by $\mnorm{\lin}$
(reserving $\norm\cdot$ for norms on
function spaces).

\subsection{Operator norms}
In general, for a linear mapping
$\lin\from{\mathbf{V}}\to\mathbf{U}$ between
normed spaces, the operator norm of $\lin$
is defined by
\begin{equation}
\mnorm\lin_o=
 \sup_{\vect}\frac{\mnorm{\lin(\vect)}}{\mnorm\vect},\quad\vect\neq0.
\end{equation}
The operator $p$-norm, $\mnorm\lin_{op}$,
$1\les{p\les\infty}$, on the space of
matrices is defined as the operator norm for
the case where the $p$-norm is used on both
$\mathbf{V}=\reals^m$ and
$\mathbf{U}=\reals^n$.  By the compactness
of the unit ball in $\reals^m$, the supremum
is attainable and
\begin{equation}
    \mnorm\lin_{op}=\max_{\mnorm\vect=1}\mnorm{\lin(\vect)}_p.
\end{equation}

 In case $\lin\from\vs\to\vs$ is a
linear transformation defined on the inner
product space $\reals^n$ equipped with the
Euclidean 2-norm, $\mnorm\lin_{o2}$ may be
calculated by
\begin{equation}
\mnorm\lin_{o2}
=\sup_{\vect,\vect'\in\vs}\frac{\abs{\lin(\vect)\cdot\vect'}}
                    {\mnorm\vect\mnorm{\vect'}},
                    \quad\vect,\vect'\neq0.
\end{equation}

The following relations hold for symmetric
matrices on $\reals^n$.
\begin{align}
    \mnorm\lin_{o1}&=\mnorm\lin_{o\infty}
     =\max_i\sum_j\abs{T_{ij}},\\
     \mnorm\lin_{o2}&=
        \max\{\abs{\lambda_1},\dots,\abs{\lambda_n}\},
\end{align}
where $\lambda_1,\dots,{\lambda_n}$ are the
(real) eigenvalues of $\lin$.  The norm
$\mnorm\cdot_{o2}$ is usually referred to as
the spectral radius norm.
\subsection{Vector norms}
We will also regard symmetric matrices as
vectors in $\reals^{{n}(n+1)/2}$ and use the
$p$-norm for them.  Thus,
\begin{equation}
    \mnorm\lin_p=\paren{\sum_{i,j}\abs{\lin_{ij}}^p}^{1/p}.
\end{equation}
In particular,
\begin{align}
\mnorm{\lin}_1
&=\sum_{i,j}\abs{\lin_{ij}},\\
\mnorm\lin_\infty&=\sup_{i,j}\abs{\lin_{ij}},\\
\mnorm\lin_2&=\sqrt{\lin_{ij}\lin_{ji}}=\paren{\sum\lambda_i^2}^{1/2}.
\end{align}
We recall that the Frobenius norm of a
matrix is
$\mnorm\lin_F=(\lin_{ij}\lin_{ij})^{1/2}$
and it is identical to $\mnorm\cdot_2$ for
symmetric matrices.
\subsection{Dual norms}
Being a finite dimensional space we may
identify the space of symmetric matrices
$\reals^{n(n+1)/2}$ with its dual space.
Thus, we may regard any symmetric matrix
$\lin$ as a linear functional so that
$\lin(S)=T_{ij}S_{ji}$ and assign to it the
dual norm $\mnorm\lin_{p^*}$
\begin{equation}
\mnorm{\lin}_{p^*}=\sup_{S}\frac{|\lin(S)|}%
    {\mnorm{S}_{p}},
\end{equation}
where we have the usual
$\mnorm\cdot_{p^*}=\mnorm\cdot_q$, for
$q=p/(p-1)$.

Dual norms may be used also for the operator
norms.  In particular, we note that
\begin{equation}
    \mnorm\lin_{o2^*}=\sum_i\abs{\lambda_i}.
\end{equation}

In closing this short review, it is noted
that the norms containing the index 2 are
associated with the Euclidean norm for
vectors and may be expressed in terms of the
eigenvalues. These norms are invariant under
orthogonal transformations of coordinates.
\subsection{The equivalence constants}
Since all the norms listed above are
equivalent, for each pair of norms
$\mnorm\cdot_a$ and ${\mnorm\cdot}_b$, there
is a finite positive number
\begin{equation}
K^a_b=\sup_\lin\frac{\mnorm\lin_a}{\mnorm\lin_b}.
\end{equation}
In particular, the following exact relations
hold:
\begin{align}
    \frac{1}{\sqrt{n}}
      &\les\frac{\mnorm\lin_{o2}}{\mnorm\lin_{o1}}
        \les\sqrt{n},\quad&
     1
      &\les\frac{\mnorm\lin_{2\phantom{o}}}{\mnorm\lin_{o2}}
        \les\sqrt{n},\quad&
        1
      &\les\frac{\hphantom{{}_{o2}}\mnorm\lin_{o2}}{\hphantom{{}_\infty}\mnorm\lin_{{\infty}}}
        \les{n},
  \\
     1
      &\les\frac{\mnorm\lin_{{o2^*}}}
        {\mnorm\lin_{{o2\phantom{{}^*}}}}
        \les{n},\quad&
       1
      &\les\frac{\mnorm\lin_{{o2^*}}}{\mnorm\lin_{2\phantom{{o}^*}}}
        \les\sqrt{n}.
\end{align}
Thus, for example, $K^{o2}_2=1$.

In the sequel, we will use $\abs{\st(x)}$ to
denote the norm of the value of a symmetric
tensor $\st$ at $x\in\gO$ using a generic
norm on the space of matrices.

\section{Optimal Boundary Conditions for
Stresses}\label{sec:OptBdCond}
Consider the following problem: given a unit
vector $\nor$ in a Euclidean 3-dimensional
space and a unit vector $\sfc$, find a
symmetric matrix $\st$ such that
\begin{enumerate}
    \item[(\textit{i})] $\st(\nor)=\sfc$---the compatibility condition;
    \item[(\textit{ii})] $\abs\st=\inf\{\abs{T},
    T(\nor)=\sfc, T=T^T\}$, i.e., $\st$ is the
    optimal symmetric matrix that satisfies
    condition (\textit{i}).
\end{enumerate}

The problem has an obvious mechanical
interpretation.  If $\nor$ denotes the
normal to the boundary at some given point,
and $\sfc$ denotes the value of the surface
traction field at that point, then,
$\st(\nor)=\sfc$ is the boundary condition
for the stress field $\st$.  Thus, a matrix
$\st$ satisfying the conditions above is the
optimal stress matrix that will satisfy the
boundary condition.  Obviously, the
normalization condition on $\sfc$ causes no
loss of generality.

Let
$\sfc_n=(\sfc\cdot\nor)\nor=\nor\otimes\nor(\sfc)$
be the normal component of $\sfc$ and let
$\sfc_t=\nor\times(\sfc\times\nor)$ be the
tangent component of $\sfc$.  Thus, denoting
the angle between $\nor$ and $\sfc$ by
$\theta$, $\abs{\sfc_n}=\cos\theta$ and
$\abs{\sfc_t}=\sin\theta$. We choose a basis
$\{f_j\}$ where $f_1=\nor$, $f_2$ is a unit
vector in the direction of $\sfc_t$ and
$f_3$ completes the other two to form a
right-hand oriented orthonormal basis.  In
this basis, the matrix of $\st$ that
satisfies the condition $\st(\nor)=\sfc$ has
to satisfy $\st_{11}=\cos\theta$,
$\st_{12}=\st_{21}=\sin\theta$, and
$\st_{13}=\st_{31}=0$. The rest of the
components cannot be determined by the
compatibility condition above and should be
determined by the requirement for minimal
norm of $\st$. (In the case where $\nor$ and
$\sfc$ are parallel, one can take any
orthonormal basis containing $\nor$.)

\subsection{Optimal boundary conditions
relative to the $\mnorm\cdot_\infty$-norm}
We wish to regard $\st$ as an element of the
dual space of symmetric matrices.  Then,
using the basis $f_i$ as above, the
compatibility condition implies that we have
a linear functional $\st_0$ defined on the
subspace $\mathbf{V}$ of symmetric matrices
containing elements of the form
\begin{equation}
[\eps]=\begin{bmatrix}
  \eps_{11} & \eps_{12} & 0 \\
  \eps_{12} & 0 & 0 \\
  0 & 0 & 0 \\
\end{bmatrix}.
\end{equation}
The functional $\st_0$ acts on elements of
$\mathbf{V}$ by
\begin{equation}
\st_0([\eps])=\eps_{11}\cos\theta+2\eps_{12}\sin\theta.
\end{equation}
Thus, on this subspace
\begin{equation}
\sup\{\abs{\st_0(\eps)}; \eps\in\mathbf{V},\
\mnorm{\eps}_1=1\}=\max\{\abs{\cos\theta},\abs{\sin\theta}\}.
\end{equation}
The extension
\begin{equation}\label{eq:OptStrMatrExpr}
[\st]=\begin{bmatrix}
  \cos\theta & \sin\theta & 0 \\
  \sin\theta & 0 & 0 \\
  0 & 0 & 0 \\
\end{bmatrix}
\end{equation}
of $\st_0$ to the space of all symmetric
matrices, has the same norm and as such it
provides the optimal boundary condition.  It
is noted that while the development depends
on the basis chosen, the optimal norm
depends only on the angle $\theta$, an
invariant quantity.

\subsection{Optimal boundary conditions
relative to the $\mnorm\cdot_2$-norm}
If one uses $\abs\st=\mnorm\st_2$ induced by
the inner product in the space of symmetric
matrices as in the previous section, the
optimal matrix can be obtained by
orthogonality conditions of the optimal
stress to basis vectors for the matrices
that correspond to the undetermined
components of the stress. This implies that
all the undetermined components should
vanish.  Thus, the optimal stress is given
in the $\{f_i\}$-basis by
Equation~\eqref{eq:OptStrMatrExpr} also,
and
\begin{equation}
    \mnorm\st_2=\sqrt{1+\sin^2\theta}.
\end{equation}

While the optimal $\st$ for the
$\mnorm\cdot_\infty$-norm depended
 on the basis $\{f_i\}$ chosen, the
 construction here is rotation-invariant.

 \subsection{The case where $\sfc=\vbase_k$}
 \label{subs:ExamOptStForEk}
 We now consider the special case where
 $\sfc=\vbase_k$, where $\vbase_k$, $k=1,2,3$ is
 a base vector.  In this case,
 $\cos\theta=\nor\cdot\vbase_k=\nor_k$ and
 $\sin\theta=\sqrt{1-\nor_k^2}$.  Using
 \begin{equation}
\st=\cos\theta\nor\otimes\nor+
 \sin\theta(\nor\otimes{f_2}+f_2\otimes\nor),
 \end{equation}
 we may write the matrix for $\st$ relative
 to the $\{\vbase_i\}$-basis.  We first note
 that
 \begin{equation}
f_2=\frac{\vbase_k-\cos\theta\nor}{\sin\theta}
   =\frac{\vbase_k-\nor_k\nor}{\sqrt{1-\nor_k^2}}\,,
 \end{equation}
 thus,

\begin{equation}
\st=\nor_k\nor\otimes\nor+
 \nor\otimes(\vbase_k-\nor_k\nor)
                +(\vbase_k-\nor_k\nor)\otimes\nor.
\end{equation}
Rearranging the terms we conclude that for
the case $\sfc=\vbase_k$, the optimal
boundary conditions for the stress are
\begin{equation}
\st=-\nor_k\nor\otimes\nor+
       (\nor\otimes\vbase_k+\vbase_k\otimes\nor),
\end{equation}
and the optimal values are
\begin{align*}
\mnorm\st_2&=\sqrt{2-\nor_k^2},\\
  \mnorm\st_\infty&=\max\braces{\abs{\nor_k},\sqrt{1-\nor_k^2}},
    \text{\ for the natural basis\
    }\{f_i\}.
\end{align*}

\subsection{Example: the 2-dimensional case for
 the $\mnorm\cdot_{o2}$-spectral radius norm}
Using the same notation as above and using
the coordinate system in the plane where the
$x$ and $y$ axes are along $f_1$ and $f_2$,
respectively, we are looking for a
$2\times2$ symmetric matrix that will
satisfy the condition $\sfc=\st(\nor)$ of
least spectral radius, i.e., minimizes
$\max\{\abs{\lambda_1},\abs{\lambda_2}\}$.
The condition $\st(\nor)=\sfc$ implies that
$\st_{xx}=\cos\theta$ and
$\st_{xy}=\sin\theta$.  Thus, to determine
$\st$ completely, one has to determine the
single number $\st_{yy}$.

In two dimensions we have the explicit
expression for the eigenvalues as
\begin{equation}
\lambda_{1,2}=\frac{\st_{xx}+\st_{yy}}2
        \pm\sqrt{\paren{\frac{\st_{xx}-\st_{yy}}2}^2
                       +{{\st_{xy}}}^2
                       }.
\end{equation}
Setting
\begin{equation}
a=\frac{\cos\theta}2=\frac{\st_{xx}}2, \quad
  z=\frac{\st_{yy}}2, \quad\text{so}  \quad
  \st_{xy}^2=\sin^2\theta=4-a^2,
\end{equation}
we have
\begin{equation}
\lambda_{1,2}=a+z\pm\sqrt{1-3a^2-2az+z^2}
 =a+z\pm\sqrt{(z-a)^2+1-4a^2}.
\end{equation}
Minimizing $\abs{\lambda_{1,2}}$ is like
minimizing $\lambda_{1,2}^2$ so we
differentiate with respect to $z$. Thus,
\begin{align*}
\frac{d\lambda_{1,2}^2}{dz}
  &=2\lambda_{1,2}\frac{d\lambda_{1,2}}{dz}\\
  &=2\lambda_{1,2}\paren{1\pm
     \frac{z-a}{\sqrt{(z-a)^2+1-4a^2}}}
\end{align*}
that vanishes identically only if $4a^2=1$.
Hence, for extremum, $\theta=0$.  In
general,
\begin{equation}
1-4a^2=1-\cos^2\theta>0,
\end{equation}
which implies that
\begin{equation}
\frac{d\lambda_{1,2}^2}{dz}
\end{equation}
has the same sign as the eigenvalue.  It
follows that $z\ne0$ can only make
$\abs{\lambda_{1,2}}$ larger and the infimum
is attained for $z=\st_{yy}=0$.

\begin{remk}
We note that the problem of optimal boundary
condition for the stress is a generalization
of the requirement in
Definition~(\ref{def:NormalVf}(\textit{i}))
for the boundary value of a normal vector
field with the difference that now we
consider matrices rather then vectors.
Indeed, $\nvf=\nor$ gives the smallest value
for $\mnorm\nvf_2$, the Euclidean norm of
the vector field $\nvf$, such that
$\nvf\cdot\nor=\nvf(\nor)=1$.
\end{remk}

\subsection{Worst case optimal boundary conditions}
By \emph{worst case optimal boundary
conditions} we refer to the inclination of
$\sfc$ relative to the normal $\nor$ for
which the norm of the optimal stress attains
a maximal value.  That is, we are looking
for
\begin{equation}
D_\|=\sup_\sfc\braces
{\inf_\st{\abs\st:\st(\nor)=t},
                       \abs\sfc=1}.
\end{equation}
The number $D_\|$ (the subscript $\|$
indicating the particular norm chosen) and
the corresponding $\st$ and $\sfc$ depend
only on the choice of norms. For example,
for the $\mnorm\cdot_2$-norm, $D_2=\sqrt2$
is attained for any unit traction $\sfc$
perpendicular to $\nor$. For the
$\mnorm\cdot_{\infty}$-norm, relative to the
boundary natural basis, $D_\infty=1$ is
attained for either traction that is
parallel to $\nor$ or traction that is
perpendicular to it.

\subsection{The spaces $\nvfs^k(\gO)$ and
                 $\vbase_k$-optimal tensor fields}
\label{subs:kOptStFields}
For $k=1,2,3$ we now choose a symmetric
tensor field $T^k=(T_{ij}^k)$ on $\bdry\gO$
that satisfies the following conditions.
\begin{enumerate}
    \item[(\textit{i})] If $\vbase_k$ denotes the $k$-th base
    vector in $\rthree$, then,
    $T^k(\nor)=\vbase_k$ so
    $T^k_{ij}\nor_j=\delta^k_i$.
    \item[(\textit{ii})]  Clearly, at each point $y$ on the boundary
    there is a collection
    of matrices that satisfy condition (\textit{i})
    above.  We choose $T^k(y)$ to be a
    symmetric matrix that satisfies
    condition (\textit{i}) and that has the least
    norm (of our choice on the space of
    matrices).  Thus,
    \begin{equation}\label{eq:OptStBC}
    |T^k(y)|=\inf_{S^k}\braces{|S^k|:S^k(\nor(y))
           =\vbase_k},
    \end{equation}
    where the infimum is taken over all
    symmetric matrices.  Thus, in the
    terminology of the preceding subsections,
    $T^k$ is the optimal stress for the
    $\vbase_k$ as boundary traction and
    the discussion of Subsection
    (\ref{subs:ExamOptStForEk}) applies.
\end{enumerate}


Clearly, the fields $T^k$ depend on $\gO$
only.  They depend continuously on $\nor$ so
by the assumption that $\bdry\gO$ is $\cee2$
they are continuous.  Consider now the worst
value of $\abs{T^k(y)}$ on the boundary,
i.e.,
\begin{equation}
\sup_{y\in\bdry\gO}\braces{|{T^k(y)}|}
 =\sup_{y\in\bdry\gO}
   \braces{\inf_{S^k}\braces{|S^k|:S^k(\nor(y))
           =\vbase_k}}.
\end{equation}
As $\bdry\gO$ is assumed to be smooth, any
angle between $\nor$ and any fixed vector is
attained on the boundary, hence, the worst
case optimal boundary conditions are
attained on the boundary always.  Thus,
\begin{equation}
\sup_{y\in\bdry\gO}\braces{|{T^k(y)}|}
 =\sup_{y\in\bdry\gO}
   \braces{\inf_{S^k}\braces{|S^k|:S^k(\nor(y))
           =\vbase_k}}=D_\|
\end{equation}
and depends only on the choice of norm.

We use the notation $T^k$ for a tensor field
on the boundary satisfying the two
conditions above. Using Whitney's extension
theorem (cf. \cite[Theorem~3.6.2]{Ziemer}),
$T^k$ can be extended to symmetric
differentiable tensor fields $\st^k$ on
$\cloo$ that satisfy the following condition
\begin{equation}\label{eq:WorstCaseAttained}
\sup_{\pis\in\cloo}\braces{|{\st^k(\pis)}|}
     =
     \sup_{y\in\bdry\gO}\braces{|{T^k(y)}|}
     =D_\|.
\end{equation}

We denote the class of symmetric
$\cee1$-tensor fields $\st^k$ on $\cloo$
that satisfy the boundary condition
$\st^k(y)=T^k(y)$, $y\in\bdry\gO$, and
condition \eqref{eq:WorstCaseAttained} above
by $\nvfs^k(\gO)$. We will refer to a tensor
field $\st^k\in\nvfs^k\ofO$ as an
$\vbase_k$-\emph{optimal tensor field}.

It is quite clear that the foregoing
discussion applies in the continuum
mechanics context to the optimal stress
field for the normalized boundary traction
$\sfc=\vbase_k$.  If $\st^k\in\nvfs^k\ofO$,
we have for the stress concentration factor
and optimal stress concentration factor,
\begin{equation}
    K_{\vbase_k,\st^k}=K_{\vbase_k}=D_\|.
\end{equation}
Although to total force on the body is not
equilibrated, $D_\|$, may serve as a bound.

\section{Bounds of the $\ldo$-Trace Operator}
\begin{thm}\label{thm:BestConstLD}
Let the constants $A\ofO$ and $B\ofO$ be
given by
\begin{align}
    A\ofO&=3D_\|, \label{eq:AforStresses}\\
    B\ofO&=
    \sum_k\norm{\nabla\cdot\st^k_0}_\infbdo,
         \label{eq:BforStresses}
\end{align}
where $\st^k_0$ is the solution of the
Dirichlet problem $\Delta\st^k_{ij}=0$, in
$\gO$, $\st^k_{ij}=T^k_{ij}$, on $\bndo$.
Then,
\begin{equation}
\norm\vf_\onebdo
 \les
                   A\ofO\norm{\eps(\vf)}_1
      +B\ofO\norm{\vf}_1,
\end{equation}
for all $\vf\in\ldo$, is an exact estimate.
\end{thm}
The following subsections present the proof.
(We will use $A$ and $B$ for $A\ofO$ and
$B\ofO$, respectively, in order to simplify
the notation.)
\subsection{The principle of virtual work}
Let $\gO$ be an open region in $\rthree$
having a smooth boundary, $\st=(\st_{ij})$ a
symmetric smooth tensor field on $\gO$ and
$\vf=(\vf_i)$ an $\ld$-vector field on
$\gO$. Then,
\begin{equation}
\st_{ij}\vf_{i,j}=(\st_{ij}\vf_i)_{,j}
   -\st_{ij,j}\vf_i.
\end{equation}
Also, by the symmetry of $\st$,
%
\begin{equation}
    \st_{ij}\vf_{i,j}
     =\st_{ij}\eps_{ij}.
\end{equation}
Thus, we may write
\begin{equation}
\ino\gO\st_{ij}\eps_{ij}\dV=\ino\gO(\st_{ij}\vf_i)_{,j}\dV
                    +\ino\gO\st_{ij,j}\vf_i\dV,
\end{equation}
and using the Green-Gauss theorem on the
first term on the right-hand side we obtain
\begin{equation}\label{eq:PrincVirtWork}
\ino\gO\st_{ij}\eps_{ij}\dV=\ino\gO\st_{ij}\half(\vf_{i,j}
                                            +\vf_{j,i})\dV
      =\ino{\bdry\gO}\st_{ij}\vf_i\nor_j\dA
                    +\ino\gO\st_{ij,j}\vf_i\dV,
\end{equation}
where $\nor$ is the unit outward pointing
normal.  We will refer to the identity above
as the \emph{principle of virtual work}.
\subsection{The bounds}
We now write the principle of virtual work
\eqref{eq:PrincVirtWork} for the vector
field $\absv\vf=(\abs{\vf_i})$. (We reserve
the notation $\abs\vect$ for the norm of the
vector $\vect$ in a finite dimensional
space.)  Thus,
\begin{equation}
\ino\gO\st_{ij}\half(|\vf_i|_{,j}+|\vf_j|_{,i})\dV
      =\ino{\bdry\gO}\st_{ij}|\vf_i|\nor_j\dA
                    +\ino\gO\st_{ij,j}|\vf_i|\dV.
\end{equation}

Let $\st^k$ satisfy $\st^k(\nor)=\vbase_k$
on $\partial\gO$ so
$\st_{ij}^k\nor_j=\delta^k_i$. Then, the
identity above assumes the form
\begin{equation}
\ino{\bdry\gO}\abs{\vf_k}\dA=
       \ino\gO\half(\abs{\vf_i}_{,j}+\abs{\vf_j}_{,i})\st^k_{ij}\dV
       -\ino\gO\st^k_{ij,j}\abs{\vf_i}\dV.
\end{equation}
Consider the integrand
\begin{equation}
\text{integrand}=\half(\abs{\vf_i}_{,j}+\abs{\vf_j}_{,i})\st^k_{ij}
\end{equation}
 of the first integral on the right.  As
 the expression is invariant under
 orthogonal transformations, it may be
 evaluated in the principle coordinate
 system of the matrix
 $(\abs{\vf_i}_{,j}+\abs{\vf_j}_{,i})/2$
 where the off-diagonal elements vanish.
 Thus, without loss of generality, we may
 write (we do not use the summation
 convention here)
 \begin{align*}
\text{integrand}&=\sum_i\abs{\vf_i}_{,i}\st^k_{ii}\\
        &
        \les\max_i\{|\st^k_{ii}|\}\sum_i\abs{\vf_i}_{,i}\\
   &\les\abs{\st^k}_{_\infty}\sum_i\abs{\vf_{i,i}}
         \quad(\text{using\ }
           \abs{\vf_i}_{,j}=\mathrm{sign}(\vf_i)\vf_{i,j})
   \\& \les\abs{\st^k}_{_\infty}\abs{\eps(\vf)}_1,
 \end{align*}
where the equality is clearly attainable.
 Thus,
\begin{align*}
\ino{\bdry\gO}\abs{\vf_k}\dA
 &\les
       \ino\gO\abs{\st^k}_{_\infty}\abs{\eps(\vf)}_1\dV
       +\ino\gO{|\st^k_{ij,j}|}\abs{\vf_i}\dV
   \\&\les
       \sup_{\pis\in\cloo}\braces{\abs{\st^k(\pis)}_{_\infty}}\ino\gO
                    \abs{\eps(\vf)}_1\dV
      +\sup_{i,x\in\cloo}\{|\st^k_{ij,j}(x)|\}\ino\gO\sum_i\abs{\vf_i}\dV.
  \end{align*}
As
\begin{equation}
    \sup_{i,x\in\cloo}\abs{\st^k_{ij,j}(x)}
    =\norm{\nabla\cdot\st^k}_\infclo,
\end{equation}
we have
\begin{equation}
\norm{\vf_k}_\onebdo=\ino{\bdry\gO}\abs{\vf_k}\dA
 \les
       \norm{\st^k}_{\infclo}
            \norm{\eps(\vf)}_1
      +\norm{\nabla\cdot\st^k}_\infclo\norm{\vf}_1,
\end{equation}
where we use
$\norm{T}_\infty=\norm{\abs{T}_\infty}_\linf$,
and $\norm{T}_1=\norm{\abs{T}_1}_\lone$, for
the respective norms of a tensor field $T$.

Adding this equation for $k=1,2,3$ we obtain
for the $\lone$-norm of the restriction of
$\vf$ to the boundary the following bound
\begin{equation}
\norm\vf_\onebdo
 \les
       \sum_k \norm{\st^k}_{\infclo}
            \norm{\eps(\vf)}_1
      +\sum_k\|\nabla\cdot\st^k\|_\infclo\norm{\vf}_1.
\end{equation}

Clearly, for
\begin{equation}
A=\inf_{\st^k}\braces{\sum_k\norm{\st^k}_\infclo:\
\st^k(\nor)=\vbase_k}
\end{equation}
and
\begin{equation}\label{eq:DefBforTensors}
B=\inf_{\st^k}\braces{\sum_k\norm{\nabla\cdot\st^k}_\infclo:
\st^k(\nor)=\vbase_k}.
\end{equation}
this bound is the tightest and we have
\begin{equation}
\norm\vf_\onebdo
 \les
                   A\norm{\eps(\vf)}_1
      +B\norm{\vf}_1, \quad\text{for all}\quad
      \vf\in\ldo.
\end{equation}

Since in general
\begin{align*}
    \sup_{\pis\in\cloo}\abs{\st^k(\pis)}_\infty
      &\ges
       \sup_{y\in\bdry\gO}\abs{\st^k(y)}_\infty
      \\&\ges
       \sup_{y\in\bdry\gO}\abs{T^k(y)}_\infty,\quad
          T^k \text{\ optimal as in
          (\ref{eq:OptStBC}})
      \\& =D_\|
\end{align*}
and equality holds for $\st^k\in
\nvfs^k(\gO)$, we conclude that $A$ is
attained for fields $\st^k\in \nvfs^k(\gO)$
and $A=3D_\|$.  This proves the first part
(Equation \eqref{eq:AforStresses}) of
Theorem~\ref{thm:BestConstLD}.
%
%
%

\subsection{Estimating $B$}

The procedure we use is completely analogous
to that of Subsection~\ref{subs:EstHarmVF}
and the proofs of theorem~\ref{thm:harVF}
and Theorem~\ref{thm:ronen}. For a tensor
field $\st$ we use $\Delta\st$ for the
Laplacian $\Delta\st_{ij}=\st_{ij,ll}$. We
say that $\st$ is \emph{harmonic} if
$\Delta\st=0$.
\begin{prop}
There is a unique harmonic tensor field
$\st^k\in\nvfs^k\ofO$, i.e., $\st^k$ is
$\vbase_k$-optimal.  For the harmonic
$\vbase_k$-optimal $\st^k$, we have
\begin{equation}\label{eq:MaxPrDivHarmSt}
    \norm{\nabla\cdot\st^k}_\infclo
     =\norm{\nabla\cdot\st^k}_\infbdo.
\end{equation}
\end{prop}
\begin{proof}
Consider the Dirichlet problem
\begin{equation}\label{eq:DirichProbStr}
    \Delta\st^k=0,
    \quad\text{in}\quad\gO,\qquad
    \st^k=T^k,\quad\text{on}\quad\bndo.
\end{equation}
The existence and uniqueness are standard.
Let $\st^k$ be harmonic, then, for each
component $\st^k_{ij}$ (no sum on repeated
indices)
\begin{equation}
\begin{split}
    \Delta((\st^k_{ij})^2)
      &=(\st^k_{ij}\st^k_{ij})_{,ll},
    \\&=
      2(\st^k_{ij,l}\st^k_{ij})_{,l}
    \\&=
      2\st^k_{ij,l}\st^k_{ij,l}+2\st^k_{ij,ll}\st^k_{ij}
    \\&\ges0.
\end{split}
\end{equation}
Thus, $\st^k_{ij}$ is subharmonic in $\gO$.
By the maximum principle for subharmonic
functions
\begin{equation}
    \max_{x\in\cloo}{(\st^k_{ij}(x))^2}
     =\max_{y\in\bndo}{(\st^k_{ij}(y))^2}
\end{equation}
and so the analogous property holds for
$\abs{\st^k_{ij}}$.  Thus, using the
boundary conditions
\begin{equation}
    \norm{\st^k}_\infclo
    =\max_{i,j,x\in\cloo}\abs{\st^k_{ij}(x)}
    =\max_{i,j,y\in\bndo}\abs{\st^k_{ij}(y)}
    =\max_{i,j,y\in\bndo}\abs{T^k_{ij}(y)}=D_\infty.
\end{equation}
Hence, the solution is also in
$\nvfs^k(\gO)$.  Finally,
\begin{equation}\begin{split}
    \Delta(\nabla\cdot\st^k)
      &=(\st^k_{ij,j})_{,ll}
    \\&=(\st^k_{ij,ll})_{,j}
    \\&=0,
\end{split}\end{equation}
so by the maximum principle for the
components of $\nabla\cdot\st^k$,
(\ref{eq:MaxPrDivHarmSt}) holds.
\end{proof}
\begin{proof}[Proof of the second part
 of Theorem~\ref{thm:BestConstLD}
 (Equation~(\ref{eq:BforStresses}))]
Clearly, as the three $\st^k$ fields are
independent, we should look for
\begin{equation}
    B_k=\inf_{\st^k}\braces{\norm{\nabla\cdot\st^k}_\infclo:
\st^k(\nor)=\vbase_k}.
\end{equation}
Thus, let $(\st^k_m)$, $m\in\mathbb{N}$, be
a sequence of $\vbase_k$-optimal tensor
fields such that
\begin{equation}
    \lim_{m\to\infty}\norm{\nabla\cdot\st^k_m}=B_k.
\end{equation}
The sets $V_m$ and $\gO_m$ may be
constructed just as in the proof of
Theorem~\ref{thm:ronen}.  Let $\st^k_{0m}$
be the solution of the Dirichlet problem in
$\gO_m$ with the boundary condition
$\st^k_{0m}(x)=\st^k_m(x)$, for
$x\in\bndo_m$, and define the harmonic
lifting accordingly as
\begin{equation}
    \cl\st^k_m(x)=\begin{cases}
    \st^k_{0m}(x),\quad&\text{for } x\in\gO_m,\\
    \st^k_{m}(x),\quad&\text{for } x\in{\cl V}_m.
    \end{cases}
\end{equation}
By the maximum principle
\begin{equation}
    \norm{\cl\st^k_{mij}}_\infclo
          \les\norm{\st^k_{mij}}_\infclo
  \quad\text{and}\quad
     \norm{\nabla\cdot\cl\st^k_{m}}_\infclo
           \les\norm{\nabla\cdot\st^k_{m}}_\infclo,
\end{equation}
so
\begin{equation}
    \lim_{m\to\infty}\norm{\nabla\cdot\cl\st^k_m}_\infclo
    =B_k.
\end{equation}
We apply the normal family argument and
uniqueness of solution as in the proof of
Theorem \ref{thm:ronen}, to obtain
\begin{equation}
    B_k=\norm{\nabla\cdot\st^k_0}_\infbdo,
\end{equation}
where $\st^k_0$ is the solution of the
Dirichlet problem \eqref{eq:DirichProbStr}.
Equation \eqref{eq:BforStresses} now follows
from \eqref{eq:DefBforTensors}.
\end{proof}
\section{Discussion}
\changed{%
We may describe the foregoing analysis in
the mechanical context.  For a given
boundary traction field, we constructed the
optimal boundary condition for the stress
field.  In Section~\ref{sec:OptBdCond}, the
case where the traction vector was a unit
vector was considered, but as the relation
between stress and traction is linear, this
causes no loss of generality.  Thus, one can
assign the optimal boundary condition for
the stress field for any given boundary
traction field.

Next, one can solve the Laplace equation for
each of the stress components.  Unlike the
usual case of continuum mechanics, we have a
unique solution without imposing
constitutive relations and the equilibrium
equations are not satisfied.  For the
harmonic solution of the boundary value
problem, the maximal stresses occur on the
boundary and these stresses are the smallest
that satisfy the traction boundary
conditions. Equation~\ref{eq:BforStresses}
and its proof indicate that the maximal
value of $\nabla\cdot\st$ is the smallest
possible.  In light of the usual equilibrium
equations $\nabla\cdot\st+\bfc=0$ of
continuum mechanics ($\bfc$ being the body
force field), we can interpret the field
$-\nabla\cdot\st$ as additional body forces
one has to supply for the equilibrium
condition to hold.  Thus, for equilibrium,
the harmonic stress field is associated with
an additional body force field whose maximum
is the least (and is attained on the
boundary). It is noted that the total of the
traction field, $\int_\bndo\sfc$, was not
required to vanish so it is not possible for
equilibrium to hold.  With the foregoing
limitations in mind, the harmonic stress
field solves the problem of optimal stress
field for a given traction.

Next, we note that the generalized stress
concentration factor may be described as the
largest optimal stress concentration factor
when we can vary the boundary traction
fields while keeping their maximal value on
the boundary to $\norm\sfc_\infbdo=1$. Thus,
in the analysis all the components of the
traction fields are set to be 1 everywhere.
Again, this precludes equilibrium as the
total force on the body in each direction is
equal to the area of the boundary.  The
mathematical analog of this limitation is
that we obtain bounds on $\trace$ and not
$\trace\qrigs$.  However, the bound on
$\norm\trace$ gives a bound
$\norm{\trace\qrigs}$ because
$\norm{\trace\qrigs}\les\norm\trace$
(Equation \eqref{eq:TraceQrigsBounded}).

}%



%
%
%
%
%
%
\begin{ack}
  The research leading to this paper was
  partially supported by the Paul Ivanier
  Center for Robotics Research and Production
  Management at Ben-Gurion University.
\end{ack}


\end{document}